# PRODUCT-LIMIT ESTIMATORS OF THE SURVIVAL FUNCTION WITH TWICE CENSORED DATA

By Valentin Patilea and Jean-Marie Rolin[1]

*CREST–ENSAI and Université Catholique de Louvain*

A model for competing (resp. complementary) risks survival data where the failure time can be left (resp. right) censored is proposed. Product-limit estimators for the survival functions of the individual risks are derived. We deduce the strong convergence of our estimators on the whole real half-line without any additional assumptions and their asymptotic normality under conditions concerning only the observed distribution. When the observations are generated according to the double censoring model introduced by Turnbull, the product-limit estimators represent upper and lower bounds for Turnbull's estimator.

**1. Introduction.** Consider the problem of nonparametric inference with competing risks survival data. The novelty we propose is that the failure time can be left-censored, for instance, at the time the study starts. For simplicity, we consider two distinct competing risks of failure, the extension to more than two competing risks being straightforward. Let $T$ and $V_1$ denote the latent independent lifetimes for each cause of failure. The failure time is $\min(T, V_1)$ and it can be censored from the left by a censoring time $U_1$. The observations are independent copies of a lifetime $Y$, a finite nonnegative random variable and a discrete random variable $A$ with values in $\{0, 1, 2\}$, where 2 indicates a left-censored failure time, while 0 and 1 correspond to an observation equal to $T$ and $V_1$, respectively. If $T$ is the lifetime of interest, we say that $Y$ is a *twice censored* observation of $T$. Associated with the problem of competing risks is the dual problem of complementary risks where the observed failure time is the maximum of the lifetimes for each cause of failure (e.g., [1]). The extension we consider here is that the

Received August 2001; revised June 2005.
[1]Supported in part by PARC Grant no. 98/03-217 of the Belgian Government.

failure time can be right-censored, for instance, at the time the experience ends.

By the plug-in (or substitution) principle applied for the empirical distribution, the nonparametric estimation of the distribution of a latent lifetime of interest is straightforward as soon as this distribution can be expressed as an explicit function of the distribution of the observed variables. The two models we propose in this paper allow for explicit inversion formulae, that is, the latent distributions of interest are explicit functionals of the distribution of the observations.

In Section 2 we introduce our latent models, while in Section 3 we provide the inversion formulae. In Section 4 we compare our model with the doubly censored data latent model proposed by Turnbull [14]. We show that the inversion formulae provide lower and upper bounds for the distribution of interest identified by Turnbull's model. Applying the inversion formulae to the empirical distribution, we deduce in Section 5 the product-limit estimators. In Sections 6 and 7 we deduce the almost sure uniform convergence and the asymptotic normality for our functionals.

**2. Latent variables models.** The random variables we consider take values in $\overline{\mathbb{R}}^+ = [0, \infty]$ endowed with $\overline{\mathcal{B}}^+$ the Borel $\sigma$-field. If $X$ is such a variable, $F_X$ denotes its distribution.

For the first latent model considered (call it Model $I$), let $T$ and $V_1$ be two lifetimes and let $U_1$ be a left-censoring time. Assume that $T$, $V_1$ and $U_1$ are independent. Suppose that $Y$ and $A$ are observed, where $Y = \max[\min(T, V_1), U_1]$ and

$$A = \begin{cases} 0, & \text{if } U_1 < T \leq V_1, \\ 1, & \text{if } U_1 < V_1 < T, \\ 2, & \text{if } \min(T, V_1) \leq U_1. \end{cases}$$

Define the observed subdistributions of $Y$ as

$$H_k(B) = P[Y \in B, A = k], \qquad k = 0, 1, 2,$$

where $B$ is a Borel set in $[0, \infty]$; the distribution of $Y$ is $H = H_0 + H_1 + H_2$. In Model $I$, the subdistributions of $Y$ can be expressed in terms of the distributions of the latent variables as follows:

(1)
$$\begin{aligned} H_0(dt) &= F_{U_1}([0,t)) F_{V_1}([t,\infty]) F_T(dt), \\ H_1(dt) &= F_{U_1}([0,t)) F_T((t,\infty]) F_{V_1}(dt), \\ H_2(dt) &= \{1 - F_T((t,\infty]) F_{V_1}((t,\infty])\} F_{U_1}(dt) \end{aligned}$$

[necessarily $H_0(\{0\}) = H_1(\{0\}) = 0$]. If $S_1 = \min(T, V_1)$ and $H_{01} = H_0 + H_1$, the three equations imply

(2) $\quad H_{01}(dt) = F_{U_1}([0,t)) F_{S_1}(dt), \qquad H_2(dt) = F_{S_1}([0,t]) F_{U_1}(dt).$

This indicates that the problem of inverting the model, that is, expressing the distributions of the latent variables in terms of the subdistributions of $Y$, can be solved in two steps. First, determine the distributions of $U_1$ and $S_1$ as in an independent left-censoring model. Next, use these distributions and the first equation in (1) to recover the distribution of $T$.

As an application of Model $I$, consider a reliability system which consists of three components $U_1$, $T$ and $V_1$, with $T$ and $V_1$ in series and $U_1$ in parallel with this series system (see, e.g., [8], Chapter 15). The lifetimes of $U_1$, $T$ and $V_1$ are independent and when the system fails we are able to determine which component failed at the same time as the system. Morales, Pardo and Quesada [9] propose the application of this model to study a certain cause of death for trees on a farm.

For our second latent model (call it Model $II$), let $U_2$ and $T$ be two lifetimes and let $V_2$ be a right-censoring time. Suppose $T$, $U_2$ and $V_2$ are independent. The observed variables are $Y$ and $A$, where $Y = \min[\max(T, U_2), V_2]$ and
$$A = \begin{cases} 0, & \text{if } U_2 < T \leq V_2, \\ 1, & \text{if } V_2 < \max(U_2, T), \\ 2, & \text{if } T \leq U_2 \leq V_2. \end{cases}$$
In Model $II$, the relationship between the subdistributions of $Y$ and the distributions of the latent variables is described by the equations

(3)
$$\begin{aligned} H_0(dt) &= F_{U_2}([0,t)) F_{V_2}([t,\infty]) F_T(dt), \\ H_1(dt) &= \{1 - F_T([0,t]) F_{U_2}([0,t])\} F_{V_2}(dt), \\ H_2(dt) &= F_T([0,t]) F_{V_2}([t,\infty]) F_{U_2}(dt) \end{aligned}$$

[necessarily $H_0(\{0\}) = 0$]. If $S_2 = \max(U_2, T)$ and $H_{02} = H_0 + H_2$, we obtain

(4) $\quad H_{02}(dt) = F_{V_2}([t,\infty]) F_{S_2}(dt), \qquad H_1(dt) = F_{S_2}((t,\infty]) F_{V_2}(dt).$

These relations show that Model $II$ can be inverted in two steps. First, as in an independent right-censoring model, recover the distributions of $V_2$ and $S_2$ from $H_{02}$ and $H_1$. Second, use the distributions of $V_2$ and $S_2$ and the first equation in (3) to determine the distribution of $T$.

Model $II$ can be interpreted as follows: consider a system consisting of three components $U_2$, $T$ and $V_2$ with independent lifetimes. Put $T$ and $U_2$ in parallel and $V_2$ in series with this parallel system (see also [2], page 767). Again, assume that we are able to determine which component failed at the same time as the system.

**3. Inversion formulae.** Recall that if $F$ is a probability distribution on $(\overline{\mathbb{R}}^+, \overline{\mathcal{B}}^+)$, the associated hazard measure is $L([0,t]) = -\ln F((t,\infty])$. Two more hazard measures can be defined,
$$L^-(dt) = \frac{F(dt)}{F([t,\infty])} \quad \text{and} \quad L^+(dt) = \frac{F(dt)}{F((t,\infty])},$$

which we call the predictable and the unpredictable hazard measure, respectively. The three hazard measures have the same continuous parts. Moreover, their point masses are in bijection: $L(\{t\}) = -\ln[1 - L^-(\{t\})] = \ln[1 + L^+(\{t\})]$. The probability distribution $F$ can be expressed as

$$F((t,\infty]) = \exp\{-L([0,t])\} = \prod_{[0,t]} (1 - L^-(ds)) = \left[\prod_{[0,t]} (1 + L^+(ds))\right]^{-1},$$

where $\prod$ is the product-integral (e.g., [6]). The mass of $L$ at infinity is irrelevant for $F$ and $F(\{\infty\}) = \exp\{-L([0,\infty))\}$.

Similarly, by reversing time, the reverse hazard measure associated to $F$ is $M((t,\infty]) = -\ln F([0,t])$. Moreover, the predictable and unpredictable reverse hazard measures are defined as

$$M^-(dt) = \frac{F(dt)}{F([0,t])} \quad \text{and} \quad M^+(dt) = \frac{F(dt)}{F([0,t))},$$

respectively. The three reverse hazard measures have the same continuous parts and their point masses satisfy $M(\{t\}) = -\ln[1 - M^-(\{t\})] = \ln[1 + M^+(\{t\})]$. We have

$$F([0,t]) = \exp\{-M((t,\infty])\} = \prod_{(t,\infty]} (1 - M^-(ds)) = \left[\prod_{(t,\infty]} (1 + M^+(ds))\right]^{-1}.$$

The mass $M(\{0\})$ is irrelevant for $F$. Moreover, $F(\{0\}) = \exp\{-M((0,\infty])\}$.

Given a nonnegative measure on $(\overline{\mathbb{R}}^+, \overline{\mathcal{B}}^+)$, we can always define a probability distribution on the same space by considering this measure as being one of $L$, $L^-$ or $L^+$ (resp. $M$, $M^-$ or $M^+$) and using the relations above. For instance, in the independent right-censoring model, one defines $L^-(dt) = H_0(dt)/H([t,\infty])$, with $H_0$ the subdistribution of the uncensored data. Then, by the equations of the model, the distribution corresponding to this $L^-$ is nothing else than the distribution of the lifetime of interest. The reverse hazard measures $M$, $M^-$ and $M^+$ are the counterparts of $L$, $L^-$ and $L^+$ to be used in left-censoring models.

We can invert our models using the hazard measures above. Since, apart from mild conditions at the origin, the inversion formulae below apply to *any* subdistributions $(H_0, H_1, H_2)$, we deduce them without any reference to the latent variables.

For inverting Model I, assume $H_0(\{0\}) = H_1(\{0\}) = 0$. In view of (2), proceed as for inverting a left-censoring model and define the predictable reverse hazard measures

$$(5) \qquad M_2^-(dt) = \frac{H_2(dt)}{H([0,t])}, \qquad M_{01}^-(dt) = \frac{H_{01}(dt)}{H([0,t)) + H_{01}(\{t\})}$$

and let $F_2^I$ and $F_{01}^I$ be the corresponding distributions. By this definition, we have $H([0,t]) = F_2^I([0,t])F_{01}^I([0,t])$. In the second step of the inversion, note that the first equation in (1) and the definition of $S_1$ imply $H_0(dt)/F_{U_1}([0,t))F_{S_1}([t,\infty]) = F_T(dt)/F_T([t,\infty])$. This suggests defining the predictable hazard measure

$$
(6) \qquad L_T^{I-}(dt) = \frac{H_0(dt)}{F_2^I([0,t))F_{01}^I([t,\infty])}.
$$

Let $F_T^I$ be its associated distribution.

For Model $II$, assume $H_0(\{0\}) = 0$. Look at the relation (4) and, exactly as in a right-censoring model, define the predictable hazard measures

$$
L_{02}^-(dt) = \frac{H_{02}(dt)}{H([t,\infty])}, \qquad L_1^-(dt) = \frac{H_1(dt)}{H((t,\infty]) + H_1(\{t\})}.
$$

Let $F_{02}^{II}$ and $F_1^{II}$ denote the corresponding distributions. Clearly, $H((t,\infty]) = F_1^{II}((t,\infty])F_{02}^{II}((t,\infty])$. In the second step of the inversion, by the first equation in (3) and the definition of $S_2$, $H_0(dt)/\{F_{V_2}([t,\infty])F_{S_2}([0,t)) + H_0(\{t\})\} = F_T(dt)/F_T([0,t])$. Consequently, define the predictable reverse hazard measure

$$
(7) \qquad M_T^{II-}(dt) = \frac{H_0(dt)}{F_1^{II}([t,\infty])F_{02}^{II}([0,t)) + H_0(\{t\})}
$$

and let $F_T^{II}$ be its associated distribution.

Now, consider the identification problem. If Model $I$ is correct, we look for conditions ensuring that $F_T^I = F_T$ on $\overline{\mathbb{R}}^+$. Define the support of $\mu$, a nonnegative measure on $[0,\infty]$, as $supp(\mu) = \{t : \mu([0,t])\mu([t,\infty]) > 0\}$. Let $B_1 = \{t : F_{U_1}([0,t))F_{V_1}([t,\infty]) > 0\}$. Deduce from (6) that the support of $L_T^{I-}$ is equal to the support of $H_0$. As $supp(H_0) = B_1 \cap supp(F_T)$,

$$F_T^I = F_T \qquad \text{on } \overline{\mathbb{R}}^+ \quad \iff \quad supp(F_T) \subset B_1.$$

By similar arguments, if Model $II$ is correct, $F_T^{II} = F_T$ on $\overline{\mathbb{R}}^+$ if and only if $supp(F_T) \subset \{t : F_{U_2}([0,t))F_{V_2}([t,\infty]) > 0\}$.

**4. Comparisons with the doubly censored data model.** The models we propose are closely related to the model for doubly (left and right) censored observations introduced by Turnbull [14]. In Turnbull's model the lifetime $T$ is independent of the censoring variables $(L,R)$ and $L \leq R$. The observations are independent copies of $Y$ and $A$, where

$$Y = \max[\min(T,R), L] = \min[\max(T,L), R],$$

$$A = \begin{cases} 0, & \text{if } L < T \leq R \text{ (no censoring)}, \\ 1, & \text{if } (L \leq)R < T \text{ (right censoring)}, \\ 2, & \text{if } T \leq L(\leq R) \text{ (left censoring)}. \end{cases}$$

If $H_k(dt) = P(Y \in dt, A = k)$, $k = 0, 1, 2$, the equations of the model are

(8) $$\begin{aligned} H_0(dt) &= \{F_L([0,t)) - F_R([0,t))\} F_T(dt), \\ H_1(dt) &= F_T((t,\infty]) F_R(dt), \\ H_2(dt) &= F_T([0,t]) F_L(dt). \end{aligned}$$

Note that the assumptions of the model imply

(9) $$H([0,t]) = F_L([0,t]) F_T([0,t]) + F_R([0,t]) F_T((t,\infty]).$$

In Turnbull's model $T$ is censored from the left by $L$ and from the right by $R$ and the observation $Y$ is always the variable in the middle. This is different from the censoring mechanisms we consider: in Model $I$ the variable $\min(T, V_1)$ is left-censored, while in Model $II$ the variable $\max(U_2, T)$ is right-censored.

Turnbull [14] proposed a nonparametric maximum likelihood estimator that can be obtained as the implicit solution of the equations (8). The implicit definition of Turnbull's estimator makes its asymptotic properties quite difficult (see [7]). Moreover, a numerical algorithm is needed for the applications.

We are interested in the relationship between our $F_T^I$, $F_T^{II}$ and $F_T$ identified by Turnbull's model. In fact, for any subdistributions $H_0$, $H_1$ and $H_2$ with $H_0(\{0\}) = H_1(\{0\}) = 0$,

$$F_T^I([0,t]) \leq F_T([0,t]) \leq F_T^{II}([0,t]) \qquad \forall\, t \geq 0,$$

where $F_T$ is the distribution of $T$ identified by Turnbull's model. Indeed, in Model $I$ use definition (6) and $H([0,t]) = F_2^I([0,t]) F_{01}^I([0,t])$ to write

$$L_T^{I-}(dt) = \frac{H_0(dt)}{F_2^I([0,t)) - H([0,t))}.$$

In Turnbull's model [relations (8) and (9)] we have

$$L_T^{-}(dt) = \frac{H_0(dt)}{F_L([0,t)) - H([0,t))}.$$

Next, the definition of $M_2^-$, the last equation in (8) and equation (9) imply

$$M_2^-(dt) = \frac{F_L([0,t]) F_T([0,t])}{F_L([0,t]) F_T([0,t]) + F_R([0,t]) F_T((t,\infty])} M_L^-(dt).$$

Deduce that the measure $M_2^-$ is smaller than the measure $M_L^-$. Therefore, $F_2^I([0,t)) \geq F_L([0,t))$, $\forall\, t \geq 0$. Hence, the measure $L_T^{I-}$ is smaller than the measure $L_T^-$, which implies $F_T^I([0,t]) \leq F_T([0,t])$, $\forall\, t \geq 0$.

On the other hand, for Model $II$, use the general relationship between $M^+$ and $M^-$, the definition (7) and $H((t,\infty]) = F_1^{II}((t,\infty]) F_{02}^{II}((t,\infty])$ and write

$$M_T^{II+}(dt) = \frac{H_0(dt)}{F_1^{II}([t,\infty]) - H([t,\infty])}.$$

Meanwhile, in Turnbull's model,

$$M_T^+(dt) = \frac{H_0(dt)}{F_R([t,\infty]) - H([t,\infty])}.$$

Next, use the definition of $L_1^-$, the general relationship between $L^+$ and $L^-$, the second equation in (3) and the equality $H((t,\infty]) = F_L((t,\infty])F_T([0,t]) + F_R((t,\infty])F_T((t,\infty])$ [this is a consequence of (9)] to deduce

$$L_1^+(dt) = \frac{F_R((t,\infty])F_T((t,\infty])}{F_L((t,\infty])F_T([0,t]) + F_R((t,\infty])F_T((t,\infty])} L_R^+(dt).$$

Clearly, the measure $L_1^+$ is smaller than the measure $L_R^+$ and, therefore, $F_1^{II}([t,\infty]) \geq F_R([t,\infty])$, $\forall t \geq 0$. Hence, the measure $M_T^{II+}$ is smaller than the measure $M_T^+$ and this implies $F_T^{II}([0,t]) \geq F_T([0,t])$, $\forall t \geq 0$.

**5. Product-limit estimators.** If we replace in the expressions of $F_T^I$ and $F_T^{II}$ the subdistributions $H_0$, $H_1$ and $H_2$ by their empirical counterparts, we obtain the product-limit estimators $F_{nT}^I$ and $F_{nT}^{II}$, respectively. For this, denote by $\{Z_j : 1 \leq j \leq M\}$ the distinct values in increasing order of $Y_i$ in a set of independent identically distributed (i.i.d.) observations $\{(Y_i, A_i) : 1 \leq i \leq n\}$. Define

$$D_{kj} = \sum_{1 \leq i \leq n} \mathbb{1}_{\{Y_i = Z_j, A_i = k\}}, \qquad N_j = \sum_{1 \leq i \leq n} \mathbb{1}_{\{Y_i \leq Z_j\}}, \qquad \overline{N}_j = \sum_{1 \leq i \leq n} \mathbb{1}_{\{Y_i \geq Z_j\}},$$

$k = 0, 1, 2$. With these definitions, the product-limit estimator of $F_T$ in Model I is

$$F_{nT}^I((Z_j, \infty]) = \prod_{1 \leq k \leq j} \left\{ 1 - \frac{D_{0k}}{U_{k-1} - N_{k-1}} \right\},$$

where

$$U_{j-1} = n \prod_{j \leq k \leq M} \left\{ 1 - \frac{D_{2k}}{N_k} \right\}.$$

The product-limit estimator of $F_T$ in Model II is given by

$$F_{nT}^{II}([0, Z_j]) = \prod_{j < k \leq M} \left\{ 1 - \frac{D_{0k}}{V_k - \overline{N}_k + D_{0k}} \right\},$$

where

$$V_j = n \prod_{1 \leq k \leq j} \left\{ 1 - \frac{D_{1k}}{\overline{N}_{k+1} + D_{1k}} \right\}.$$

When the doubly censored data model is considered, our product-limit estimators represent lower and upper bounds for Turnbull's estimator. These bounds may serve for the numerical algorithms used to compute Turnbull's estimator.

**6. Strong convergence.** We study the strong (almost sure or a.s.) uniform convergence of $F_{nT}^I$ and $F_{nT}^{II}$. Since, in fact, the estimators $F_{nT}^I$ and $F_{nT}^{II}$ are built as explicit functionals of the empirical distribution, we deduce their asymptotic behavior, in particular, the strong convergence, whatever the properties of the underlying censoring mechanism are. Hereafter, we use the following rule: the subscript $n$ indicates the empirical version of the quantities we consider. Moreover, if $\mu$ is a nonnegative measure on $(\overline{\mathbb{R}}^+, \overline{\mathcal{B}}^+)$ and $f$ is a measurable function, $\mu(f) = \int f(t)\mu(dt)$.

For the strong convergence, we recall a result of Rolin [12], an extension of the strong law under right-censorship proved by Stute and Wang [13]. Let $H = \sum_{1 \leq r \leq g} H_r$ be a probability distribution decomposed into $g$ subdistributions. If $\mathcal{I} \subset \mathcal{K} = \{r : 1 \leq r \leq g\}$, let $H_\mathcal{I} = \sum_{r \in \mathcal{I}} H_r$. For $\mathcal{J}_k \subset \mathcal{I}_k \subset \mathcal{K}$, $k = 1, 2$, define

$$L_k^-(dt) = \frac{H_{\mathcal{J}_k}(dt)}{H((t, \infty]) + H_{\mathcal{I}_k}(\{t\})}$$

and consider the measure $G(dt) = \exp\{-L_2([0, t))\} L_1^-(dt)$.

THEOREM 6.1. *If $G(f) < \infty$, then $G_n(f) \to G(f)$ a.s. and in the mean.*

The same result holds if we define the predictable reverse hazard measures

$$M_k^-(dt) = \frac{H_{\mathcal{J}_k}(dt)}{H([0, t)) + H_{\mathcal{I}_k}(\{t\})}, \qquad k = 1, 2,$$

and consider the measure $G(dt) = \exp\{-M_2((t, \infty])\} M_1^-(dt)$.

Let us extend the number of hazard measures associated with Model $I$ by defining

$$M_0^-(dt) = \frac{H_0(dt)}{H([0, t)) + H_{01}(\{t\})}, \qquad M_1^-(dt) = \frac{H_1(dt)}{H([0, t)) + H_1(\{t\})}.$$

Consider $F_0^I, F_1^I$, the corresponding distributions. Deduce that $M_{01} = M_0 + M_1$, where $M_0, M_1$ and $M_{01}$ are the reverse hazard measures associated with $M_0^-, M_1^-$ and $M_{01}^-$ [see (5)], respectively. In view of equations (2), deduce $H([0, t)) = F_2^I([0, t)) F_{01}^I([0, t))$ and $H_{01}(\{t\}) = F_2^I([0, t)) F_{01}^I(\{t\})$. Therefore,

$$H_0(dt) = F_2^I([0, t)) F_{01}^I([0, t]) M_0^-(dt).$$

Consequently, in the expression of the predictable hazard measure defining $F_T^I$ [see (6)], we get rid of $F_2^I$ and obtain

$$L_T^{I-}(dt) = \frac{F_{01}^I([0, t])}{F_{01}^I([t, \infty])} M_0^-(dt).$$

THEOREM 6.2. *If $f$ is a nonnegative Borel measurable function defined on $(\overline{\mathbb{R}}^+, \overline{\mathcal{B}}^+)$ such that $L_T^{I-}(f) < \infty$, then, almost surely as $n \to \infty$, $L_{nT}^{I-}(f) \to L_T^{I-}(f)$.*

Theorem 6.2 is a direct consequence of the following lemma.

LEMMA 6.3. (i) *If $L_T^{I-}(f\mathbb{1}_{[0,t]}) < \infty$ and $F_{01}^I([t,\infty]) > 0$, then a.s.*

$$L_{nT}^{I-}(f\mathbb{1}_{[0,t]}) \to L_T^{I-}(f\mathbb{1}_{[0,t]}), \qquad n \to \infty.$$

(ii) *If $L_T^{I-}(f\mathbb{1}_{[t,\infty]}) < \infty$ and $F_2^I([0,t)) > 0$, then almost surely as $n \to \infty$,*

$$L_{nT}^{I-}(f\mathbb{1}_{[t,\infty]}) \to L_T^{I-}(f\mathbb{1}_{[t,\infty]}).$$

PROOF. (i) First, Theorem 6.1 implies that any empirical distribution function defined by the empirical reverse hazard measures of Model $I$ converges uniformly on $[0,\infty]$. Now,

$$\left| L_{nT}^{I-}(f\mathbb{1}_{[0,t]}) - \int_{(0,t]} \frac{f(s)}{F_{01}^I([s,\infty])} F_{n01}^I([0,s]) M_{n0}^-(ds) \right|$$

$$\leq \frac{\|F_{n01}^I - F_{01}^I\|}{F_{n01}^I([t,\infty])} \int_{(0,t]} \frac{f(s)}{F_{01}^I([s,\infty])} F_{n01}^I([0,s]) M_{n0}^-(ds).$$

The second member of the inequality tends to zero almost surely because $F_{n01}^I([t,\infty]) \to F_{01}^I([t,\infty]) > 0$ a.s. and, by Theorem 6.1 applied for $G(ds) = \exp\{-M_{01}((s,\infty])\}M_0^-(ds)$,

$$\int_{(0,t]} \frac{f(s)}{F_{01}^I([s,\infty])} F_{n01}^I([0,s]) M_{n0}^-(ds) \to L_T^{I-}(f\mathbb{1}_{[0,t]}), \qquad \text{a.s.}$$

(ii) First, looking at the definition of $M_{01}^-$, by a simple computation,

$$H_{01}([s,\infty]) \leq F_{01}^I([s,\infty]) \leq \frac{H_{01}([s,\infty])}{H_{01}([0,\infty])}.$$

Using definition (6) for the predictable hazard measure defining $F_T^I$, we have

$$\left| L_{nT}^{I-}(f\mathbb{1}_{[t,\infty]}) - \int_{[t,\infty]} \frac{f(s)}{F_2^I([0,s))} \frac{H_{n0}(ds)}{F_{n01}^I([s,\infty])} \right|$$

$$\leq \frac{\|F_{n2}^I - F_2^I\|}{F_{n2}^I([0,t))} \int_{[t,\infty]} \frac{f(s)}{F_2^I([0,s))} \frac{H_{n0}(ds)}{F_{n01}^I([s,\infty])}$$

$$\leq \frac{\|F_{n2}^I - F_2^I\|}{F_{n2}^I([0,t))} \int_{[t,\infty]} \frac{f(s)}{F_2^I([0,s))} \frac{H_{n0}(ds)}{H_{n01}([s,\infty])}.$$

Now, almost surely $F_{n2}^I([0,t)) \to F_2^I([0,t))$, which is strictly positive. Since

$$\int_{[t,\infty]} \frac{f(s)}{F_2^I([0,s))} \frac{H_0(ds)}{H_{01}([s,\infty])} \leq H_{01}([0,\infty])^{-1} L_T^{I-}(f \mathbb{1}_{[t,\infty]}) < \infty,$$

a new application of Theorem 6.1 provides the result. $\square$

Denote by $t_{0k}$ the left endpoint and by $t_{1k}$ the right endpoint of the support of $H_k$, $k = 0, 1, 2$. We have the following corollary of Theorem 6.2. Note that the strong uniform convergence of $F_{nT}^I$ is obtained without any additional assumption, apart from that of i.i.d. observations and the condition $H_0(\{0\}) = H_1(\{0\}) = 0$.

COROLLARY 6.4. (a) If $L_T^{I-}([0, t_{10})) < \infty$, then, almost surely,

$$\sup_{0 \leq t < t_{10}} |L_{nT}^I([0,t]) - L_T^I([0,t])| \to 0$$

and $L_{nT}^I(\{t_{10}\}) \to L_T^I(\{t_{10}\})$. If $L_T^{I-}([0,t_{10})) = \infty$, then, almost surely,

$$\sup_{0 \leq s \leq t} |L_{nT}^I([0,s]) - L_T^I([0,s])| \to 0$$

for all $t < t_{10}$ and $L_{nT}^I([0, t_{10})) \to \infty$.
(b) Almost surely, $\|F_{nT}^I - F_T^I\| = \sup_{0 \leq t \leq \infty} |F_{nT}^I([0,t]) - F_T^I([0,t])| \to 0$.

PROOF. The Glivenko–Cantelli theorem provides the result in (a) with $L_T^I$ and $L_{nT}^I$ replaced by $L_T^{I-}$ and $L_{nT}^{I-}$, respectively. The similar result for the hazard measure $L_{nT}^I$ is obtained by taking care of the fact that $L_T^I(\{t_{10}\}) = \infty$ if $L_T^{I-}(\{t_{10}\}) = 1$. This happens if $t_{10} \geq t_{11}$, $H_0(\{t_{10}\}) > 0$ and $H_1(\{t_{10}\}) = 0$. The convergence of $F_{nT}^I$ is implied by the convergence of the associated hazard measure $L_{nT}^I$. $\square$

The strong uniform convergence of $F_{nT}^{II}$ can be obtained in a similar way. Define

$$L_0^-(dt) = \frac{H_0(dt)}{H([t,\infty])}, \qquad L_2^-(dt) = \frac{H_2(dt)}{H((t,\infty]) + H_1(\{t\}) + H_2(\{t\})}$$

and consider $F_0^{II}, F_2^{II}$, the corresponding distributions. After some manipulations we can get rid of $F_1^{II}$ in the definition (7):

$$M_T^{II-}(dt) = \frac{F_2^{II}([t,\infty]) F_0^{II}([t,\infty])}{1 - F_2^{II}([t,\infty]) F_0^{II}((t,\infty])} L_0^-(dt).$$

Next, apply Theorem 6.1 (see [10] for the details).

**7. Asymptotic normality.** Let $(D[a,b], \|\cdot\|)$ be the space of càdlàg functions defined on $[a,b] \subset [0,\infty]$, endowed with the supremum norm. $BV_C[a,b] \subset D[a,b]$ is the set of càdlàg functions with total variation bounded by $C$. The integrals with respect to functions which are not of bounded variation have to be understood via partial integration. Finally, weak convergence is denoted by $\rightsquigarrow$ and is in the sense considered by Pollard [11], that is, $D[a,b]$ is endowed with the ball $\sigma$-field.

Given the explicit form of $F_{nT}^I$ and $F_{nT}^{II}$, a convenient approach for proving weak convergence is the delta method (e.g., [5] and [15], Section 3.9). For proving Hadamard differentiability, the denominators appearing in the maps used to define $F_T^I$ and $F_T^{II}$ should stay away from zero. Therefore, we have to complete the delta method with a tool for treating the endpoints of the intervals on which weak convergence is finally proved.

LEMMA 7.1 ([11], page 70). *Let $X, X_1, X_2, \ldots$ be random elements of $(D[a,b], \|\cdot\|)$ with the distribution of $X$ concentrated on a separable set. Suppose, for each $\varepsilon, \delta > 0$ there exist approximating random elements $AX, AX_1, AX_2, \ldots$ such that $AX_n \rightsquigarrow AX$, $P(\|X - AX\| > \delta) < \varepsilon$ and*

$$(10) \qquad \limsup_{n \to \infty} P(\|X_n - AX_n\| > \delta) < \varepsilon.$$

*Then $X_n \rightsquigarrow X$.*

For brevity, we consider only the asymptotic normality of $F_{nT}^I$; similar arguments apply for $F_{nT}^{II}$. The empirical central limit theorem yields

$$\sqrt{n}(H_n - H, H_{0n} - H_0, H_{2n} - H_2) \rightsquigarrow (G, G_0, G_2) \qquad \text{in } D^3([0,\infty]).$$

Now, we prove that $\sqrt{n}(M_{n2}^- - M_2^-)$ and $\sqrt{n}(F_{n2}^I - F_2^I)$ converge weakly to Gaussian limits. The computation of the covariance structures for the limit processes in this section is elementary, albeit tedious (see [10] for some formulae).

LEMMA 7.2. *Let $M_{2t}^- = M_2^-((t,\infty])$ and $M_{n2t}^-$ be the corresponding estimator. Assume that*

$$(11) \qquad \int_{(t_{00},\infty]} \frac{M_2^-(du)}{H([0,u])} = \int_{(t_{00},\infty]} \frac{H_2(du)}{H([0,u])^2} < \infty,$$

*where $t_{00} = \inf\{t : H_0([0,t]) > 0\}$. Then*

$$(12) \quad \sqrt{n}(H_n - H, H_{n0} - H_0, M_{n2}^- - M_2^-) \rightsquigarrow (G, G_0, G_M) \qquad \text{in } D^3[t_{00}, \infty],$$

*where $(G, G_0, G_M)$ is a zero-mean Gaussian process with*

$$(13) \qquad G_{Mt} = \int_{(t,\infty]} \frac{dG_{2u}}{H([0,u])} - \int_{(t,\infty]} \frac{G_u}{H([0,u])^2} H_2(du).$$

Moreover, if $F_{2t}^I = F_2^I([0,t])$, then

$$\sqrt{n}(H_n - H, H_{n0} - H_0, F_{n2}^I - F_2^I) \rightsquigarrow (G, G_0, G_3) \qquad in \ D^3[t_{00}, \infty],$$

where $(G, G_0, G_3)$ is a zero-mean Gaussian process with

(14) $$G_{3t} = F_2^I([0,t]) \int_{(t,\infty]} \frac{dG_{Mu}}{1 - M_2^-(\{u\})}.$$

PROOF. The map $(A, B) \to \int_{(\cdot,\infty]} (1/A) \, dB$ is Hadamard-differentiable on a domain of the type $\{(A, B) : A \in D[a,b], B \in BV_C[a,b], A \geq \epsilon\}$, $C, \epsilon > 0$, at every point such that $1/A$ is of bounded variation. The derivative map is given by $(\alpha, \beta) \to \int_{(\cdot,\infty]} (1/A) \, d\beta - \int_{(\cdot,\infty]} (\alpha/A^2) \, dB$. Therefore, the delta method for the map $(H, H_0, H_2) \to (H, H_0, M_2^-)$ yields the weak convergence of $\sqrt{n}(H_n - H, H_{n0} - H_0, M_{n2}^- - M_2^-)$ in $D^3[\sigma, \infty]$, provided that $H([0,\sigma]) > 0$.

For the weak convergence in $D^3[t_{00}, \infty]$, consider the pathwise limit of $G_{M\sigma}$ as $\sigma \downarrow t_{00}$, which exists in view of (11). It remains to verify (10) when $H([0,t_{00}]) = 0$. It suffices to prove the following: (a) for any $\varepsilon, \delta > 0$, there exists $\sigma = \sigma(\varepsilon, \delta) > t_{00}$ such that

(15) $$\limsup_{n \to \infty} P\left( \sup_{U \leq t \leq \sigma} \sqrt{n} |M_{n2}^-([t,\sigma)) - M_2^-([t,\sigma))| > \delta \right) < \varepsilon;$$

and (b) $\sqrt{n} M_2^-((t_{00}, U)) \to 0$, in probability, where $U = \min_i Y_i$. To ensure (a), reverse the time and apply the arguments usually used to check the "tightness at $\tau_H = \sup\{t : H([0,t]) < 1\}$" when proving weak convergence for Nelson–Aalen and Kaplan–Meier estimators (see [3], Theorem 6.2.1, [4]). For (b), first note that (11) ensures $M_2^-((t_{00}, \infty])$ is finite. This implies $F_2^I([0,t_{00}]) > 0$ (use, e.g., arguments as in Lemma 6 of [6]). Since in general $M^-$ is smaller than $M$, deduce

$$M_2^-((t_{00}, U)) \leq M_2((t_{00}, U)) = \ln \frac{F_2^I([0,U))}{F_2^I([0,t_{00}])} \leq \frac{F_2^I((t_{00}, U))}{F_2^I([0,t_{00}])}.$$

Let $u_n^\lambda = \sup\{s : \sqrt{n} F_2^I((t_{00}, s)) \leq \lambda\}$ (see also [16]). We have

$$P(\sqrt{n} F_2^I((t_{00}, U)) > \lambda) \leq P(U > u_n^\lambda) = H((u_n^\lambda, \infty])^n$$
$$\leq \{1 - F_2^I((t_{00}, u_n^\lambda]) F_{01}^I([0, u_n^\lambda])\}^n$$
$$\leq \left(1 - \frac{\lambda^2}{n} \frac{F_{01}^I([0, u_n^\lambda])}{F_2^I((t_{00}, u_n^\lambda])}\right)^n \to 0.$$

The convergence to zero is true because, in view of (11),

$$\frac{F_2^I((t_{00}, u_n^\lambda])}{F_{01}^I([0, u_n^\lambda])} \leq \int_{(t_{00}, u_n^\lambda]} \frac{F_2^I(ds)}{F_{01}^I([0,s])} = \int_{(t_{00}, u_n^\lambda]} \frac{M_2^-(ds)}{H([0,s])} \to 0, \qquad n \to \infty.$$

Now, (b) is clear. For the last part of the lemma, apply the delta method for the map $A \to \pi_{(\cdot,\infty]}(1 - A(ds))$ defined on $BV_C[t_{00}, \infty]$, for some $C > 0$. □

REMARK. In view of the variance of the process $G_3$, it seems possible to relax condition (11) when $F_2^I([0, t_{00}]) = 0$ (see also [4]). However, in the following, due to the lack of an obvious martingale structure for $L_{nT}^{I-} - L_T^{I-}$, it is convenient to keep the denominator appearing in the definition of $L_T^{I-}$ away from zero when $t \downarrow t_{00}$. For this, we have to impose $F_2^I([0, t_{00}]) > 0$ and, in this case, (11) is needed to bound the variance of $G_{3t}$ when $t \downarrow t_{00}$.

Now we state the asymptotic normality for $L_{nT}^{I-}$ and $F_{nT}^I$. The notation $A_-$ means that we consider the left-limits of the process $A$.

THEOREM 7.3. *Suppose condition* (11) *holds. Let* $t_{00} < \tau$ *such that* $H_{01}([0, \tau)) < 1$. *If* $L_{Tt}^{I-} = L_T^{I-}([0, t])$, *then* $\sqrt{n}(L_{nT}^{I-} - L_T^{I-}) \rightsquigarrow V$ *in* $D[0, \tau]$, *where*

$$V_t = \int_{(0,t]} \frac{dG_{0u}}{(F_2^I - H)([0, u))} - \int_{(0,t]} \frac{G_{3u-} - G_{u-}}{(F_2^I - H)^2([0, u))} H_0(du), \qquad t \in [0, \tau],$$

*is a zero-mean Gaussian process. Moreover, if* $F_{Tt}^I = F_T^I([0, t])$, *then we have* $\sqrt{n}(F_{nT}^I - F_T^I) \rightsquigarrow W$ *in* $D[0, \tau]$, *with* $W$ *the zero-mean Gaussian process*

$$W_t = F_T^I((t, \infty]) \int_{(0,t]} \frac{dV_u}{1 - L_T^{I-}(\{u\})}.$$

PROOF. Since $F_{01}^I([\tau, \infty]) > 0$ and, by (11), $F_2^I([0, t_{00}]) > 0$, we have $\inf_{(t_{00}, \tau]}(F_2^I - H)([0, s)) > \epsilon$, for some $\epsilon > 0$. Thus, if $H_0(\{t_{00}\}) = 0$, the weak convergence of $\sqrt{n}(L_{nT}^{I-} - L_T^{I-})$ is obtained by the delta method for the map $(A, B) \to \int_{(t_{00}, \cdot]}(1/A_-) \, dB$ (see [15], pages 382–384).

When $H_0(\{t_{00}\}) > 0$ (hence, necessarily $t_{00} > 0$), in the definition of $L_T^{I-}$, we also have to take into account $F_2^I([0, t_{00}))$. For this, extend the weak convergence in (12) on $D^3[0, \infty]$ by considering a modified predictable reverse hazard function

$$M_{2t}^- = M_2^-((t, \infty]) = \int_{(t, \infty]} \frac{H_2(du)}{H([0, u \vee t_{00}])}, \qquad t \in [0, \infty].$$

Let $M_{n2}^-$ be the empirical counterpart. Since the denominator in the last display stays away from zero, the weak convergence of $\sqrt{n}(H_n - H, H_{n0} - H_0, F_{n2}^I - F_2^I)$ in $D^3[0, \infty]$ is easily obtained by the delta method, where now $F_2^I$, $F_{n2}^I$ correspond to the modified $M_2^-$, $M_{n2}^-$, respectively. Note that now $F_2^I([0, t_{00})) > 0$. The processes $G_M$ and $G_3$ are still defined according to (13) and (14), respectively. Since the modification of $M_{2t}^-$ and $M_{n2}^-$ does not change the definitions of $L_T^{I-}$ and $L_{nT}^{I-}$, the delta method yields the weak convergence of $\sqrt{n}(L_{nT}^{I-} - L_T^{I-})$. The last part of the theorem is obtained by the delta method for the product-integration map. $\square$

**Acknowledgments.** The authors wants to thank the Editor, an Associate Editor and the referees for a careful reading of the manuscript. The valuable suggestions of Stephen Lagakos are gratefully acknowledged.

CREST–ENSAI
Rue Blaise Pascal BP 37203
35172 Bruz cedex
France
E-mail: patilea@ensai.fr

Institut de Statistique
Université Catholique de Louvain
20 voie du Roman Pays
1348 Louvain-la-Neuve
Belgique
E-mail: rolin@stat.ucl.ac.be